\documentclass[11pt, onesided, reqno]{amsart}

\usepackage{graphicx}
\usepackage{amsmath,amsthm}
\usepackage{amsfonts,amssymb}
\usepackage{enumerate}
\usepackage{color}

\newtheorem{theorem}{Theorem}
\newtheorem{proposition}[theorem]{Proposition}
\newtheorem{corollary}[theorem]{Corollary}
\newtheorem{lemma}[theorem]{Lemma}

\newtheorem{THEA}{Theorem A\!\!}

\newtheorem{THEB}{Theorem B\!\!}

\newtheorem{THEC}{Theorem C\!\!}

\theoremstyle{definition}

\newtheorem{remark}[theorem]{Remark}

\newcommand{\E}{\mathbb{E}}

\newcommand{\R}{\mathbb{R}}
\newcommand{\Pb}{\mathbb{P}}

\renewcommand{\L}{\mathcal{L}}
\def\tLL{{\widetilde L}}

\newcommand{\cD}{\mathcal{D}}
\newcommand{\cN}{\mathcal{N}}

\newcommand{\M}{\mathbf{M}}
\renewcommand{\i}{{\rm i}}

\def\lacc{\left\{}
\def\lcr{\left[}
\def\lpa{\left(}
\def\lva{\left|}
\def\M{{\bf M}}
\def\racc{\right\}}
\def\rpa{\right)}
\def\rcr{\right]}
\def\rva{\right|}

\newcommand{\Marg}{\M_{\alpha,\rho,\gamma}}
\newcommand{\Margb}{\M_{\alpha,\rho,\gamma}^{(\beta)}}
\renewcommand{\a}{\alpha}

\newcommand{\Un}{{\bf 1}}

\def\={{\;\mathop{=}\limits^{\text{(law)}}\;}}

\def\elaw{\stackrel{d}{=}}

\newcommand*\pFqskip{8mu}
\catcode`,\active
\newcommand*\pFq{\begingroup
        \catcode`\,\active
        \def ,{\mskip\pFqskip\relax}%
        \dopFq
}
\catcode`\,12
\def\dopFq#1#2#3#4#5{%
        {}_{#1}F_{#2}\biggl[\genfrac..{0pt}{}{#3}{#4};#5\biggr]%
        \endgroup
}

\allowdisplaybreaks

\begin{document}

\title{Cram\'er's estimate for stable processes with power drift}

\bigskip

\author[Christophe Profeta]{Christophe Profeta}

\address{Laboratoire de Math\'ematiques et Mod\'elisation d'\'Evry, Universit\'e d'Evry-Val d'Essonne, B\^atiment IBGBI, 23 Boulevard de France, 91037 Evry Cedex, France. {\em Email}: {\tt christophe.profeta@univ-evry.fr}}

\author[Thomas Simon]{Thomas Simon}

\address{Laboratoire Paul Painlev\'e, Universit\'e de Lille,  Cit\'e Scientifique, 59655 Villeneuve d'Ascq Cedex, France. {\em Email}: {\tt simon@math.univ-lille1.fr}}



\begin{abstract} We investigate the upper tail probabilities of the all-time maximum of a stable L\'evy process with a power negative drift. The asymptotic behaviour is shown to be exponential in the spectrally negative case and polynomial otherwise, with explicit exponents and constants. Analogous results are obtained, at a less precise level, for the fractionally integrated stable L\'evy process. We also study the lower tail probabilities of the integrated stable L\'evy process in the presence of a power positive drift.
\end{abstract}

\maketitle

\section{Introduction and statement of the results}

\bigskip

Let $L$ be a real strictly $\alpha$-stable L\'evy process with characteristic exponent
$$\Psi(\lambda) \;=\;\log(\E[e^{\i \lambda L_1}])\; =\;-\; (\i \lambda)^\alpha e^{-\i\pi\alpha\rho\, {\rm sgn}(\lambda)},\qquad \lambda \in \R,$$
where $\alpha \in (0,2]$ is the scaling parameter and $\rho=\Pb[L_1>0]$ is the positivity parameter. Recall e.g. from Lemma 14.11 and Theorem 14.19 in \cite{Sato} that $\rho\in[1-1/\a,1/\a]$ for $\a\in[1,2]$ and $\rho\in[0,1]$ for $\a\in(0,1),$ and that with this normalization, for $\a\in(0,2)$ the L\'evy measure of $L$ has density  
$$\nu(x)\; =\; \frac{\Gamma(1+\alpha)}{\pi}\left(\frac{\sin(\pi\alpha (1-\rho))}{|x|^{1+\alpha}}\Un_{\{x<0\}}\; +\; \frac{\sin(\pi\alpha \rho) }{x^{1+\alpha}}\Un_{\{x>0\}}\right).$$
Throughout, we assume that $L$ takes positive values i.e. $\rho\neq 0,$ and we exclude the degenerate case $\a = \rho =1$ where $L$ is a unit drift. With these  restrictions, $L$ has positive jumps if and only if $\a \le 1$ or $\a > 1$ and $\rho <1/\a.$\\

Consider the positive random variable
$$\Marg\; =\; \sup_{t\ge 0} \{L_t - t^\gamma\}.$$
It is well-known from e.g. Proposition 48.10 in \cite{Sato} that 
$$\Pb[\Marg < \infty]\; = \; \left\{\begin{array}{ll} 1 & \mbox{if $\gamma\a > 1$}\\
0 & \mbox{if $\gamma\a \le 1.$}\end{array}\right.$$
In this paper, we are concerned with the asymptotic behaviour of  
$$\Pb[\Marg \geq x], \qquad x\rightarrow \infty,$$
in the relevant case $\gamma\a > 1.$ In the literature, the evaluation of such asymptotics having various applications in insurance is coined as Cram\'er's estimate. In the case of a linear drift $\gamma =1,$ we refer to (XI.6.16) and (XII.5.10) in \cite{Fel} for random walks and to \cite{BD} for L\'evy processes having one-sided exponential moments. Applied to stable L\'evy processes, the main result of \cite{BD} shows 
\begin{equation}
\label{One}
\Pb[\M_{\a,1/\a,1} \geq x]\;\sim \; e^{-x}
\end{equation}
for $\a > 1,$ and it is well-known that the asymptotics is in fact an equality - see \cite{Zo} or Corollary VII.2 in \cite{Ber}. For more general power drifts and a class of Gaussian processes fulfilling a certain scaling property, we refer to \cite{HP} which, applied to the important case of Brownian motion with a parabolic drift, yields
\begin{equation}
\label{Two}
\Pb[\M_{2,1/2,2} \geq x]\;\sim \; \frac{1}{\sqrt{3}} \exp\lacc - \frac{4}{3\sqrt{3}}\, x^{3/2}\racc.
\end{equation}
Let us mention that this estimate has been recently refined in Theorem 2.1 of \cite{GT}, where a complete asymptotic expansion at infinity is obtained - see also Lemma 2.1 and the references therein for closed expressions of the density of $\M_{2,1/2,2}$ in terms of the Airy function. The first result of the present paper is the following general estimate, extending (\ref{One}) and (\ref{Two}).

\begin{THEA}\label{theo:1} Assume $\gamma \alpha>1$. 

\medskip 

{\em (a)}  If $L$ has positive jumps, one has
$$\Pb[\Marg \geq x]\; \sim\; \frac{\sin(\pi\alpha \rho)}{\pi}\,\Gamma(\alpha - 1/\gamma)\Gamma(1+1/\gamma) \;x^{\frac{1}{\gamma}-\a}.$$

{\em (b)}  If $L$ does not have positive jumps, one has  
$$\Pb[\M_{\a,1/\a,\gamma} \geq x]\; \sim\;\sqrt{\frac{\a -1}{\gamma\a -1}}\, \exp\lacc-(\alpha-1)\,\gamma^{\frac{\alpha}{\alpha-1}}\,(\gamma\a -1)^{\frac{1-\gamma\a}{\gamma(\alpha-1)}}\,x^{\frac{\gamma\a-1}{\gamma(\alpha-1)}} \racc.$$
\end{THEA}

\bigskip

In the specific case $\a\in (1/2,2]$ and $\gamma =2,$ these estimates are somehow reminiscent of those previously obtained in \cite{BerBur} in the framework of Burgers turbulence with stable noise initial data. More precisely, if we set 
$$\underline{\bf M}_{\alpha,\rho,\gamma}^{[x]}\; =\; \sup_{t\in[0,x]} \{L_t - t^\gamma\}\qquad\mbox{and}\qquad \overline{\bf M}_{\alpha,\rho,\gamma}^{[x]}\; =\; \sup_{t\ge x} \{L_t - t^\gamma\},$$
then the main result of \cite{BerBur} states that
$$\Pb[\overline{\bf M}_{\alpha,\rho,2}^{[x]}\ge \underline{\bf M}_{\alpha,\rho,2}^{[x]}]\; \asymp\; x^{1-2\a}$$
if $L$ has positive jumps\footnote{Here and throughout, we use the standard notation $f(x)\asymp g(x)$ to express the fact that there exist two positive finite constants $\kappa_1, \kappa_2$  such that $\kappa_1 f(x)\leq g(x) \leq \kappa_2 f(x)$ as $x\to\infty$ or as $x\to 0$, the nature of the limit being clear from the context.}, and that 
$$\log \Pb[\overline{\bf M}_{\alpha,1/\a,2}^{[x]}\ge \underline{\bf M}_{\alpha,1/\a,2}^{[x]}]\; \sim\; -\kappa_\a\, x^{\frac{2\a-1}{\alpha-1}}$$
for some explicit $\kappa_\a\in(0,\infty)$ if $L$ does not have positive jumps.
 Roughly speaking, when $x$ is large the event
$$\lacc \overline{\bf M}_{\alpha,\rho,2}^{[x]}\ge \underline{\bf M}_{\alpha,\rho,2}^{[x]}\racc$$ 
amounts to the fact that the translated process $L_t - L_x - (t-x)^2$ crosses a level of size $x^2$ for some $t\ge x,$ which explains heuristically why the asymptotics of 
$$\Pb[\overline{\bf M}_{\alpha,\rho,2}^{[x]}\ge \underline{\bf M}_{\alpha,\rho,2}^{[x]}]\qquad\mbox{and}\qquad \Pb[\M_{\a,\rho,2} \ge x^2]$$ 
are comparable. Our arguments are quite different from those of \cite{BerBur}. They rely on the compensation formula for the case with positive jumps and on some ad hoc and rather involved estimates combined with Laplace's method in the spectrally negative case. One might wonder if such arguments could not help refine the results of \cite{BerBur}, but we have not investigated this question.\\

In the second part of the paper we consider the Riemann-Liouville (or fractionally integrated) stable process with parameter $\beta > 0$, defined as 
$$L^{(\beta)}_t\;=\; \int_0^t (t-s)^\beta dL_s\; =\; \beta\int_0^t (t-s)^{\beta-1} L_s\, ds, \qquad t\ge 0.$$
The process $\{L^{(\beta)}_t, \, t\ge 0\}$ is stable in the broad sense of \cite{SamTaq}, and by Proposition 3.4.1 therein we have
\begin{equation}
\label{Variance} L^{(\beta)}_1\;\elaw\; (1+\alpha\beta)^{-1/\alpha} L_1.
\end{equation}
Recall also that $\{L^{(\beta)}_t, \, t\ge 0\}$ is self-similar with index $\beta +1/\a,$ non-Markovian, and that it has a.s. continuous sample paths. Consider the positive random variable
$$\Margb\; =\; \sup_{t\ge 0} \{L^{(\beta)}_t - t^{\beta+\gamma}\},$$
and observe from e.g. Theorem 10.5.1 in \cite{SamTaq} and self-similarity that$$\Pb[\Margb < \infty]\; =\; \Pb[\Marg < \infty]$$ 
for every $\beta > 0.$  As a rule, the non-Markovian character of a given process makes however its passage times across a level more difficult to investigate. Our second main result has a less precise character. 

\begin{THEB}\label{theo:2} Assume $\gamma\alpha>1$.

\medskip 

{\em (a)} If $L$ has positive jumps, one has 
$$\Pb[\Margb \geq x]\; \asymp\;  x^{\frac{1-\gamma\a}{\beta +\gamma}}.$$

{\em (b)} If $L$ does not have positive jumps,  one has
$$\log\Pb[\M_{\a,1/\a,\gamma}^{(\beta)} \geq x] \;\sim\; -c_{\a,\beta,\gamma}\, x^{\frac{\gamma\alpha -1}{(\alpha-1)(\gamma+\beta)}}$$
with $c_{\a,\beta,\gamma} = (\alpha-1)\,(\gamma+\beta)^{\frac{\alpha}{\alpha-1}}\, (\a\beta+1)^{\frac{\gamma+\beta-1-\alpha\beta}{(\alpha-1)(\gamma +\beta)}}\,(\gamma\a -1)^{\frac{1-\gamma\a}{(\alpha-1)(\gamma +\beta)}}\, >\, 0.$
\end{THEB}

 The method to get these estimates differs here for the lower bound and the upper bound. The former uses a simple scaling argument, inspired by that of \cite{HP}, and amounts to a comparison with the upper tails of $L_1.$ The latter relies on telescoping sums for the case with positive jumps, and on a simple yet powerful association lemma in the spectrally negative case - see Lemma \ref{lem:ST}.\\

In the last part of the paper, we study the lower tail problem for the integrated stable process with a power positive drift. In a Gaussian framework, lower tail probabilities have many applications described in \cite{LS}. In a self-similar framework they are connected to the persistence probabilities, whose applications are also manifold - see the recent survey \cite{AS}. We show the following.

\begin{THEC}\label{theo:3} 
Assume $\gamma\alpha > 1$ and $\rho\in (0,1).$ For every $\mu\ge 0,$ one has
$$ \Pb\lcr\sup_{0\leq t\leq 1} \lacc L_t^{(1)} + \mu t^{1+\gamma}\racc\leq\, \varepsilon\rcr\; \asymp\;\varepsilon^{\frac{\alpha\rho}{(\alpha+1)(\a(1-\rho)+1)}}. 
$$
\end{THEC}

Above, we have excluded the case $\rho =1,$ where the estimate amounts by monotonicity to the one-dimensional estimate $\Pb[L_1 +\mu \leq\varepsilon],$ which is exponentially small for $\mu = 0$ - see e.g. (14.35) in \cite{Sato} - and zero for $\mu > 0.$ Theorem C is an extension of Theorem A in \cite{PSInt} which dealt with the case $\mu = 0.$ In this respect, we should mention that the condition $\gamma\a >1$ on the drift power is optimal: in the Cauchy case $\a=\gamma =1,$ the same Theorem A in \cite{PSInt} shows that the lower tail probability exponent depends on $\mu.$ Our argument relies in an essential way on the strong Markov property of the bidimensional process $\{(L^{(1)}_t, L_t), \, t\ge 0\}$ and is hence specific to the case $\beta = 1.$ The other cases are believed to be challenging. To give but one example, for $\a=\beta =2$ and $\mu = 0$, finding the right asymptotics of 
$$ \Pb\lcr\sup_{0\leq t\leq 1} \lacc L_t^{(2)}\racc\leq\, \varepsilon\rcr$$
as $\varepsilon\to 0$ is still an open problem on Brownian motion - see Section 3.3 in \cite{AS}. In our proof the aforementioned association Lemma \ref{lem:ST} plays also a significant role. Unfortunately, its one-sided character prevents us from dealing with the case of a negative power drift. We leave this question, whose connection to Burgers turbulence with stable L\'evy process initial data in the case $\a > 1$ and $\gamma =1$ is described in Section 4.1 of \cite{AS}, to future research.   

\section{Proof of Theorem A}

\subsection{The case with positive jumps} We will use the standard notation
 $$c_+\; =\;\frac{\Gamma(1+\alpha)}{\pi} \sin(\pi\rho\alpha)\; > \; 0$$
for simplicity.  Defining for every $x>0$ the stopping time
$$T_x\;=\; \inf\{t\geq 0;\, L_t> t^{\gamma} +x\},$$
we have $\Pb[\M_{\a,\rho,\gamma} \ge x]\; =\; \Pb[T_x<\infty].$ We also set $K_x=L_{T_x} - T^\gamma_x -x$ for the overshoot at $T_x.$ For every $f : \R^+\to\R^+$ measurable and such that $f(0)=0$, the compensation formula - see \cite{Ber} p. 7 or Theorem 19.2 in \cite{Sato} - implies
\begin{eqnarray*}
\E\left[f(K_x)\, \Un_{\{T_x<\infty\}}\right]&= &\E\lcr \sum_{t\geq 0}f\left(L_{t^-} + \Delta L_t -t^\gamma -x\right)  \Un_{\{L_u < u^\gamma +x \;\forall u < t, \, t^\gamma +x < L_{t^-}+\Delta L_t\}}\rcr\\
&= & c_+\,\E\lcr \int_0^{\infty} \!\! dt \int_0^{\infty}\!\! f\left(L_t + s -t^\gamma-x\right)  \Un_{\{L_u < u^\gamma +x \;\forall u < t, \, t^\gamma +x < L_{t}+s\}}\,s^{-1-\alpha} ds \rcr \\
&= & c_+\, \int_0^{\infty}\!\! dt  \int_0^{\infty}\!\!\E\lcr f(z-t^\gamma-x) \Un_{\{L_u < u^\gamma +x \;\forall u < t, \, t^\gamma +x < z\}} (z-L_t)^{-1-\alpha} \rcr dz. 
\end{eqnarray*}
Taking $f(u)=\Un_{\{u>0\}}$ and integrating in $z$, we obtain 
\begin{eqnarray*}
\Pb[K_x >0,\,  T_x<\infty]\!\! &=&\!\! \frac{c_+}{\alpha} \int_0^{\infty}  \E\left[(t^\gamma+x- L_t)^{-\alpha}   \,\Un_{\{L_u < u^\gamma +x\;\forall u < t\}}\right]\, dt \\
&=&\!\!\frac{c_+}{\alpha}\lpa \int_0^{\infty} \E\left[(s^\gamma+1-  x^{\frac{1}{\gamma\alpha}-1}L_s)^{-\alpha}   \,\Un_{\{ x^{\frac{1}{\gamma\alpha}-1}  L_v < v^\gamma +1\;\forall v <s\}}\right]  ds\rpa x^{\frac{1}{\gamma}-\alpha}\\
& \sim & \!\!\frac{c_+}{\alpha}\lpa \int_0^{\infty} (s^\gamma+1)^{-\alpha} \,ds\rpa x^{\frac{1}{\gamma}-\alpha}\\ 
& \sim &  \frac{\sin(\pi\alpha \rho)}{\pi}\,\Gamma(\alpha - 1/\gamma)\Gamma(1+1/\gamma) \;x^{\frac{1}{\gamma}-\a}
\end{eqnarray*}
where the second equality follows by scaling, the convergence on the third line is obtained by bounded and monotone convergence (decomposing into $\{L_s < 0\}$ and $\{L_s \ge 0\}$ inside the expectation), and the evaluation of the integral on the fourth line is standard. To conclude the proof, it remains to show that $L$ does not creep at $T_x$, in other words that
\begin{equation}
\label{Creep}
\Pb[K_x=0,\,T_x<\infty]\; =\; 0.
\end{equation}
The latter is in accordance with the well-known fact that $L$ does not creep at a fixed level $x >0$ - see Theorem VI.19 and Lemma VIII.1 in \cite{Ber}. However, this result does not apply here since we consider the first passage time above a moving boundary. To show (\ref{Creep}), fix $x > 0$ and decompose 
$$\Pb[L_s \geq s^\gamma +x]\;  =\; P_1(s)\, +\, P_2(s)$$
for every $s \ge 0,$ with
$$\lacc\begin{array}{lll}
P_1(s)& =& \Pb[\tLL_{s-T_x} + T_x^\gamma \geq s^\gamma, \,K_x = 0,\, T_x<s]\\
P_2(s)& =&\Pb[\tLL_{s-T_x} + L_{T_x}\geq s^\gamma +x,\, K_x > 0,\, T_x<s],\end{array}\right.$$
where $\tLL$ is a copy of $L$ which is independent of $(T_x, L_{T_x})$, by the strong Markov property. On the one hand, we see by scaling and e.g. Property 1.2.15 in \cite{SamTaq} that
$$ \Pb[L_s \geq s^\gamma +x] \;\sim\; \frac{c_+}{\alpha}\, s^{1-\gamma\alpha}.$$
On the other hand, we have
$$P_1(s) \; \geq \; \Pb[\tLL_{s-T_x} \geq s^\gamma, \,K_x = 0,\, T_x<s/2]\, \geq \,  \Pb[\tLL_1 \geq 2^{\frac{1}{\alpha}} s^{\gamma-\frac{1}{\alpha}}]\; \Pb[K_x = 0,\, T_x< s/2]$$
and passing to the limit, we obtain 
$$\liminf_{s\rightarrow \infty}  s^{\gamma\alpha -1}\,P_1(s) \; \geq\;  \frac{c_+}{2\alpha}\;\Pb[K_x = 0,\, T_x<\infty].$$
Hence, we see that (\ref{Creep}) is a consequence of
\begin{equation}
\label{Creepy}
P_2(s) \;\sim\; \frac{c_+}{\alpha}\, s^{1-\gamma\alpha}.
\end{equation}
Applying the compensation formula as above, we obtain
$$P_2(s) \;=\; c_+\int_0^s dt \int_0^{\infty}\Pb[\tLL_{s-t} +L_t+z\geq s^\gamma +x, \,L_t + z> t^\gamma +x, \, L_u < u^\gamma + x\;\forall u < t]\, z^{-1-\alpha}\, dz.$$
Changing the variables $z=s^{\gamma} y$ and $t=su,$ we see that $c_+^{-1} s^{\gamma\alpha -1}P_2(s)$ equals  
$$\int_0^1\!\! du \int_0^{\infty}\Pb[s^{\frac{1}{\alpha}} \tLL_{1-u} + s^{\frac{1}{\alpha}}L_u+s^{\gamma} y \geq s^\gamma +x, \, s^{\frac{1}{\alpha}}L_u > s^\gamma (u^\gamma-y) +x, \, s^{\frac{1}{\alpha}} L_u < s^\gamma\,u^\gamma + x\;\forall u < 1]\,y^{-1-\alpha}dy,$$
which converges as $s\to\infty$ to
$$\int_0^1 du \int_1^{\infty} y^{-1-\alpha} dy \; =\;\frac{1}{\alpha}\cdot$$
This shows (\ref{Creepy}), and completes the proof.
 
\qed

\begin{remark} (a) Setting, here and throughout, $L_t^* = \sup\{L_s, \, s\in[0,t]\}$ for every $t > 0,$ we have
$$\lim_{\gamma\rightarrow \infty} \Pb[\Marg \geq x]\; =\; \Pb[L_1^* \geq x]$$
for every $x\ge 0.$ Passing formally to the limit $\gamma\to\infty$ in Theorem A (a), we can infer
\begin{equation}
\label{Bing}
\Pb[L_1^* \geq x]  \;\sim\; \frac{\Gamma(\alpha)  \sin(\pi\alpha \rho)}{\pi}\, x^{-\a}
\end{equation}
which is a standard and rigorous estimate - see Theorem 10.5.1 in \cite{SamTaq} and Proposition VIII.4 in \cite{Ber}.\\

(b) Taking  $f(u) = \Un_{\{u\geq r x\}}$ for some $r>0$ and applying as above the compensation formula leads to the estimate 
$$\Pb[K_x\geq r x ,\, T_x<\infty] \; \sim\; \frac{c_+}{\alpha}  \lpa \int_0^{\infty}(r + u^\gamma+1)^{-\alpha}\,  du\rpa x^{\frac{1}{\gamma}-\a}\;\sim\; (r+1)^{\frac{1}{\gamma} -\a}\Pb[T_x <\infty].$$
This implies the following limit theorem for the law of the renormalized overshoot:
$$\L\lpa x^{-1} K_x\,\bigg| \,T_x<+\infty\rpa\;  \to\;   {\rm Pareto}\, (\alpha-1/\gamma) \qquad \mbox{as $x\to\infty$.}$$
This observation seems new even in the classical case of a linear drift $\gamma = 1$ with $\a > 1.$ Notice that still in the case of a linear drift, the limit behaviour of the overshoot is very different for L\'evy processes having finite exponential moments. If we consider for example the tempered stable subordinator with negative unit drift and L\'evy measure having density
$$\nu(x)\; =\; \frac{\a\,e^{-cx}}{\Gamma(1-\a) x^{\a+1}}\,\Un_{\{x > 0\}}$$
for some $c\in (0,1),$ then we are in the framework of \cite{BD} with $\omega \in (0,1)$ and $\mu^* < \infty$ so that $C > 0$ in (5) therein. By Remark 2 of \cite{BD}, this implies that $K_x$ converges at infinity to some proper random variable. \\

(c) In the case $\a > 1, \rho = 1-1/\a$ and $\gamma =1,$ the Laplace transform of ${\M}_{\a,1-1/\a,1}$ can be computed with the help of Zolotarev's well-known general formula in \cite{Zo}: one finds
$$\E[e^{-\lambda{\M}_{\a,1-1/\a,1}}]\; =\; \frac{1}{1 +\lambda^{\a-1}}\cdot$$  
This Laplace transform can be easily inverted and yields the identity in law
$${\M}_{\a,1-1/\a,1}\;\elaw\; {\bf L}^{\frac{1}{\a-1}}\,\times\, {\bf Z}_{\a-1}$$
where ${\bf L}\sim{\rm Exp}(1)$ and ${\bf Z}_{\a-1}$ has a standard positive $(\a-1)-$stable law with Laplace transform
$$\E[e^{-\lambda {\bf Z}_{\a-1}}]\; =\; e^{-\lambda^{\a-1}},$$
both random variables being independent. This shows that the law of ${\M}_{\a,1-1/\a,1}$ is the so-called Mittag-Leffler distribution of parameter $\a-1$ which is studied e.g. in Exercise 34.4 of \cite{Sato} - see also the references therein. In particular, there exists a closed expression for the survival function of ${\bf M}_{\a,1-1/\a,1}$ in terms of the classical Mittag-Leffler function, which leads to a complete and simple asymptotic expansion at infinity: one has
$$\Pb[{\M}_{\a,1-1/\a,1}\, > \,x] \; =\; E_{\a -1}(-x^{\a-1})\; \sim\;\sum_{n\ge 1} \frac{(-1)^{n-1} x^{-(\a-1) n}}{\Gamma (1- (\a-1) n)}$$
where we have used the standard expansion 18.1(7) in \cite{EMOT}. Observe from the complement formula for the Gamma function that the first term matches the one that can be derived from Theorem A (a), in this specific case. Notice also the following closed formula for the distribution function, as a convergent series:
$$\Pb[{\M}_{\a,1-1/\a,1}\, \le \,x] \; =\; 1\, -\, E_{\a -1}(-x^{\a-1})\; = \;\sum_{n\ge 1} \frac{(-1)^{n-1} x^{(\a-1)n}}{\Gamma (1+ (\a-1) n)}\cdot$$ 
Let us finally refer to \cite{Fur} for related results in the presence of a compound Poisson process.\\

\end{remark}

\subsection{The case with no positive jumps} Applying the strong Markov property at $T_x$ and using the absence of positive jumps, we get
\begin{eqnarray*}
\int_0^{\infty} (1-e^{-\lambda t})\, \Pb[L_t> t^\gamma +x]\,dt & = & \E\lcr\Un_{\{T_x<\infty\}} \int_0^{\infty} (1-e^{-\lambda (T_x+t)})\,\Un_{\{L_t+T_x^\gamma> (t+T_x)^\gamma\}}\, dt\rcr\\
& = & \E\lcr\Un_{\{T_x<\infty\}} \int_0^{\infty} (1-e^{-\lambda (T_x+t)})\,\Un_{\{t^{1/\a}L_1^+ > (t+T_x)^\gamma - T_x^\gamma\}}\, dt\rcr
\end{eqnarray*}
where we have set $a^+ = \max(a,0)$ and, on the right-hand side, $L$ and $T_x$ are independent. Integrating both sides on $(0,\infty)$ with respect to $\lambda^{-\nu-1} d\lambda$ with $\nu\in (0,1)$, we deduce 
\begin{eqnarray*}
\int_0^{\infty} t^{\nu}\,\Pb[L_t> t^\gamma +x]\,dt &= &\E\lcr\Un_{\{T_x<\infty\}} \int_0^{\infty} (T_x+t)^{\nu}\, \Un_{\{t^{1/\alpha} L_1^+ > (t+T_x)^\gamma - T_x^\gamma\}}\, dt\rcr\\
&= &\E\left[\Un_{\{T_x<\infty\}} \, T_x^{1+\nu}\, \int_0^{\infty} (1+t)^{\nu}\, \Un_{\{L_1^+ >  T_x^{\gamma-1/\alpha}\,\varphi_{\a,\gamma} (t)\}}\, dt\rcr,
\end{eqnarray*}
where
$$\varphi_{\a,\gamma}(t)\; =\;  \frac{(1+t)^\gamma -1}{t^{1/\alpha}}$$
is an homeomorphism from $(0,\infty)$ to $(0,\infty),$ because $\alpha\gamma>1$ and $\alpha>1$. This implies the identity
\begin{eqnarray}
\nonumber \int_0^{\infty}\!\! t^{\nu}\,\Pb[L_1^+> t^{-1/\a}(t^\gamma +x)]\,dt\!\! &= &\!\!\E\lcr \Un_{\{T_x<\infty\}} \,T_x^{1+\nu}  \int_0^{\infty} (1+t)^{\nu}\, \Un_{\{\varphi_{\a,\gamma}^{-1}(T_x^{1/\a - \gamma} L_1^+)>  t\}}\, dt\rcr\\
\label{eq:1}\!\!& = &\!\!\frac{1}{1+\nu}\, \E\lcr\Un_{\{T_x<\infty\}} \,T_x^{1+\nu} \lpa(1+\varphi_{\a,\gamma}^{-1}( T_x^{1/\a-\gamma}L_1^+))^{1+\nu} -1\rpa\rcr
\end{eqnarray}
which extends to all $\nu>-1$ by analyticity, since $L_1^+$ has moments of every order. We will now study the asymptotic behaviour of both sides of (\ref{eq:1}), introducing the crucial parameter
$$\nu_0\; = \; \frac{\alpha(\gamma-1)}{\alpha-1}\; >\; -1.$$
We begin with the left-hand side, which is easy. 
\begin{lemma}
\label{lem:Lap}
One has
$$\int_0^{\infty}t^{\nu}\,\Pb[L_1^+> t^{-1/\a}(t^\gamma +x)]\,dt \,\sim\, \frac{\gamma^{\frac{\a}{1-\a}} ((\gamma\a -1)x)^{\frac{\nu -\nu_0}{\gamma}}}{\sqrt{(\a -1)(\gamma\a-1)}}\,\exp\lacc -(\alpha-1)\,\gamma^{\frac{\alpha}{\alpha-1}}\,(\gamma\a -1)^{\frac{1-\gamma\a}{\gamma(\alpha-1)}}\, x^{\frac{\alpha\gamma-1}{\gamma(\alpha-1)}}\racc.$$
\end{lemma}

\proof
By (14.35) in \cite{Sato}, we have the asymptotic behaviour
$$p_1(x)\; \sim\;\frac{\alpha^{-\frac{1}{2(\alpha-1)}}}{\sqrt{2\pi(\alpha-1)}}\,x^{\frac{2-\alpha}{2(\alpha-1)}} \exp\lacc-(\alpha-1)\, \a^{\frac{\a}{1-\a}}\, x^{\frac{\alpha}{\alpha-1}}\racc$$
at infinity, where $p_1$ stands for the density of the random variable $L_1$. Making a change of variable and applying Watson's lemma - see also Theorem 2.5.3 in \cite{Z}, this easily implies
\begin{equation}
\label{eq:iGamma}\Pb[L^+_1> x] \;\sim\;\frac{\alpha^{\frac{1}{2(\alpha-1)}}}{\sqrt{2\pi(\alpha-1)}}\,   x^{-\frac{\alpha}{2(\alpha-1)}}\, \exp \lacc-(\alpha-1)\, \a^{\frac{\a}{1-\a}}\, x^{\frac{\alpha}{\alpha-1}}\racc. 
\end{equation}
On the other hand, we can rewrite
\begin{equation}
\label{eq:iii}
\int_0^{\infty}t^{\nu}\,\Pb[L_1^+> t^{-1/\a}(t^\gamma +x)]\,dt\; =\; x^{\frac{\nu+1}{\gamma}}  \int_0^{\infty}  s^\nu \,\Pb[L^+_1> x^{\frac{\alpha\gamma-1}{\alpha\gamma}}\eta(s)]\,ds\end{equation}
where $\eta(s)= s^{-1/\a}(s^\gamma+1)$ reaches its global minimum on $(0,\infty)$ at $s_\ast= (\alpha\gamma -1)^{-1/\gamma},$ with
$$\eta(s_\ast)\; =\; \gamma\a(\gamma\a-1)^{\frac{1-\gamma\a}{\gamma\a}}\qquad\mbox{and}\qquad \eta''(s_\ast)\; =\;\frac{\gamma(\gamma\a-1)^{2+1/\a}}{\a}\cdot$$
Plugging (\ref{eq:iGamma}) into the right-hand side of (\ref{eq:iii}), we obtain
\begin{multline*}\int_0^{\infty}t^{\nu}\,\Pb[L_1^+> t^{-1/\a}(t^\gamma +x)]\,dt \\\sim\;\frac{\alpha^{\frac{1}{2(\alpha-1)}}}{\sqrt{2\pi(\alpha-1)}}  \, x^{\frac{\nu+(1-\nu_0)/2}{\gamma}}     \int_0^{\infty} s^\nu\,\eta(s)^{\frac{\alpha}{2(1-\alpha)}}\, \exp\lacc-(\alpha-1)\,\a^{\frac{\alpha}{1-\alpha}}  \,\eta(s)^{\frac{\alpha}{\alpha-1}}  x^{1+\nu_0}\racc\,ds,
 \end{multline*}
which yields the required asymptotic behaviour, by  Laplace's method.

\endproof

We will now analyze the right-hand side of (\ref{eq:1}), which is more involved. Introducing the function 
$$\Phi_{\a,\gamma,\nu}(x)\; =\; x^{-\frac{(1+\nu)\alpha}{\gamma\a-1}} \left((1+\varphi_{\a,\gamma}^{-1}(x))^{1+\nu}-1\right)$$
on $(0,\infty),$ we can rewrite (\ref{eq:1}) as 
\begin{equation}
\label{eq:1bis} \int_0^{\infty} t^{\nu}\,\Pb[L_1^+> t^{-1/\a}(t^\gamma +x)]\,dt\; =\;\frac{1}{1+\nu}\, \E\lcr\Un_{\{T_x<\infty\}}\left(L_1^+\right)^{\frac{(1+\nu)\alpha}{\gamma\a-1}} \Phi_{\a,\gamma,\nu} (T_x^{1/\a-\gamma}\,L_1^+)\rcr 
\end{equation}
Taking $\nu =\nu_0$ and observing that $\varphi_{\a,\gamma}^{-1}(t)\sim (t/\gamma)^{\frac{\alpha}{\alpha-1}}$ as $t\to 0$ and $\varphi_{\a,\gamma}^{-1}(t)\sim t^{\frac{\alpha}{\gamma\a-1}}$ as $t\to\infty,$ we get 
$$\lim_{x\rightarrow 0}\Phi_{\a,\gamma, \nu_0}(x)\; =\; (1+\nu_0)\gamma^{\frac{\alpha}{1-\alpha}}\; >\; 0
\qquad \text{ and }\qquad \lim_{x\rightarrow \infty}\Phi_{\a,\gamma, \nu_0}(x)\; =\; 1.$$
Therefore, since $\Phi_{\a,\gamma,\nu_0}$ is continuous and positive on $(0,\infty),$ we have
 
$$0\;<\; \inf_{x > 0}\lacc \Phi_{\a,\gamma,\nu_0}(x)\racc \;<\; \sup_{x > 0}\lacc\Phi_{\a,\gamma,\nu_0}(x)\racc\; <\;\infty.$$
Going back to  (\ref{eq:1bis}) and using Lemma \ref{lem:Lap}, we finally get the crude asymptotics

\begin{equation}
\label{eq:leqleq}
\Pb[T_x<\infty]\;\asymp\;\exp\lacc -(\alpha-1)\,\gamma^{\frac{\alpha}{\alpha-1}}\,(\gamma\a -1)^{\frac{1-\gamma\a}{\gamma(\alpha-1)}}\, x^{\frac{\alpha\gamma-1}{\gamma(\alpha-1)}}\racc.
\end{equation}
In order to obtain an exact asymptotics and finish the proof, we will need the following technical lemma.
\begin{lemma}\label{lem:Phi}
For every $\nu \in (-1, \gamma-1/\a -1],$ the function $\Phi_{\a,\gamma,\nu}$ is an homemorphism from $(0,\infty)$ to $(0,1).$
\end{lemma}
\proof
First, it is easy to see from the aforementioned asymptotics of $\varphi_{\a,\gamma}$ at zero and infinity that
$$\lim_{x\rightarrow 0}\Phi_{\a,\gamma, \nu}(x)\; =\; 0
\qquad \text{ and }\qquad \lim_{x\rightarrow \infty}\Phi_{\a,\gamma, \nu}(x)\; =\; 1$$
for $\nu \in (-1, \gamma-1/\a -1],$ and it is plain that $\Phi_{\a,\gamma, \nu}$ is continuous. Since $\varphi_{\a,\gamma}$ increases on $(0,\infty),$ we are reduced to show that
$$z\;\mapsto\;\Phi_{\a,\gamma,\nu}\left(\varphi_{\a,\gamma}(z)\right)\; = \;\frac{\left((1+z)^{1+\nu}-1\right)\,z^{\frac{(1+\nu)}{\gamma\a-1}} }{\left((1+z)^\gamma-1\right)^{\frac{(1+\nu)\alpha}{\gamma\a-1}}}$$
increases on $(0,\infty).$ Setting $y=(1+z)^{\gamma}-1$ and $f_c(x) =  (1+x)^c - x^c,$ we obtain 
$$\left(\Phi_{\a,\gamma,\nu}\left(\varphi_{\a,\gamma}(z)\right) \right)^{\frac{\gamma\a-1}{1+\nu}}\; =\;\left(f_{\frac{1+\nu}{\gamma}}(y^{-1})  \right)^{\frac{\gamma\a-1}{1+\nu}} f_{\frac{1}{\gamma}}(y^{-1})$$
which, since $f_c$ decreases for $c\in(0,1],$ shows that $\Phi_{\a,\gamma, \nu}$ increases for $\gamma\geq 1$ and $\nu \in (-1, \gamma-1]$. Assuming last $\gamma<1,$ we need to prove that 
$$x\;\mapsto\; g_{\a,\gamma,\nu}(x)\; =\;  \left(f_{\frac{1+\nu}{\gamma}}\left(x\right)  \right)^{\alpha-\frac{1}{\gamma}} \left(f_{\frac{1}{\gamma}}\left(x\right)\right)^{\frac{1+\nu}{\gamma}}$$
decreases on $(0,\infty)$. Setting $c= \frac{1+\nu}{\gamma} \in (0,1),$ we compute
$$g_{\a,\gamma,\nu}^\prime(x)\; =\; \frac{c\, g_{\a,\gamma,\nu}(x)}{\gamma(1+x)}\left(\alpha \gamma - (\alpha \gamma-1)  \frac{x^{c-1}}{f_{c}(x)} -  \frac{x^{\frac{1}{\gamma}-1}}{f_{\frac{1}{\gamma}}(x)}\right)\; < \; \frac{c\, g_{\a,\gamma,\nu}(x)}{\gamma(1+x)}\left(\alpha \gamma - (\alpha \gamma-1)  \frac{x^{c-1}}{f_{c}(x)}\right).$$
It is easy to see that $x\mapsto x^{1-c}f_{c}(x)$ increases from $(0,+\infty)$ to $(0,c),$ and we finally obtain
$$g^\prime_{\a,\gamma,\nu}(x) \; <\; \frac{((c-1)\gamma\a +1)g_{\a,\gamma,\nu}(x)}{\gamma(1+x)}\; \le \; 0$$
as soon as $\nu\le \gamma-1/\alpha-1$.

\endproof

\begin{corollary}
\label{eq:A}
For every $A\ge 0,$ one has
$$\frac{\Pb[T_x\leq A]}{\Pb[T_x<+\infty]} \; \to \; 0$$
as $x\to\infty.$
\end{corollary}

\proof
Set $\nu=\varepsilon - 1$ with $\varepsilon>0$ small enough for $\Phi_{\a,\gamma,\varepsilon -1}$ to increase on $(0,\infty).$ By (\ref{eq:1bis}) and the fact that $L_1^+$ and $T_x$ are independent, we have
\begin{eqnarray*}
\varepsilon \int_0^{\infty} t^{\varepsilon - 1}\,\Pb[L_1^+> t^{-1/\a}(t^\gamma +x)]\,dt & = & \E\lcr\Un_{\{T_x<\infty\}}\left(L_1^+\right)^{\frac{\alpha}{\gamma\a-1}} \Phi_{\a,\gamma,\varepsilon-1} (T_x^{1/\a-\gamma}\,L_1^+)\rcr \\
&\geq & \E\left[\left(L_1^+\right)^{\frac{\alpha}{\gamma\a-1}}  \Phi_{\a,\gamma,\varepsilon-1} (A^{1/\a-\gamma}\,L_1^+)  \right]\,  \Pb[ T_x\leq A].
\end{eqnarray*}
Combining now the crude asymptotics (\ref{eq:leqleq}) and Lemma \ref{lem:Lap}, we deduce that there exists $K>0$ such that
$$ \frac{\Pb[T_x\leq A]}{\Pb[T_x<+\infty]}\; \leq \; K\, x^{\frac{\varepsilon- 1 -\nu_0}{\gamma}}\;\to\; 0$$
as $x\to\infty,$ taking $\varepsilon > 0$ small enough. 

\endproof

We can now finish the proof. Taking $\nu =\nu_0$ in (\ref{eq:1bis}), we first decompose the quantity
$$\gamma^{\frac{\alpha}{\a-1}}\int_0^{\infty}t^{\nu_0}\,\Pb[L_1^+> t^{-1/\a}(t^\gamma +x)]\,dt$$
into
$$\E[(L_1^+)^{\frac{\alpha}{\alpha-1}}]\,\Pb[T_x<\infty]\; +\; \frac{1}{\Phi_{\a,\gamma,\nu_0}(0+)}\,\E\lcr  \Un_{\{T_x<\infty\}} (L_1^+)^{\frac{\alpha}{\alpha-1}} \lpa\Phi_{\a,\gamma,\nu_0}( T_x^{1/\a -\gamma}L_1^+) -\Phi_{\a,\gamma,\nu_0}(0+)\rpa \rcr.$$
Applying Lemma \ref{lem:Lap} and the moment evaluation
$$\E[(L_1^+)^{\frac{\alpha}{\alpha-1}}]\; =\; \frac{1}{\a-1}$$
which is e.g. a consequence of (2.6.20) in \cite{Z}, we see that the proof will be complete as soon as 
\begin{equation}
\label{llinf}
\E\lcr  \Un_{\{T_x<\infty\}} (L_1^+)^{\frac{\alpha}{\alpha-1}} \lpa\Phi_{\a,\gamma,\nu_0}( T_x^{1/\a -\gamma}L_1^+) -\Phi_{\a,\gamma,\nu_0}(0+)\rpa \rcr\; \ll\; \Pb[T_x <\infty], \quad x\to\infty.
\end{equation}
But, decomposing according to $\{T_x\leq A\}$ or $\{T_x> A\}$, the left-hand side is bounded by 
\begin{multline*}
\frac{2}{\a-1}\, \sup_{z > 0}\lacc\Phi_{\a,\gamma,\nu_0}(z)\racc\,\Pb[T_x\leq A]\\ 
+\;  \E\left[(L_1^+)^{\frac{\alpha}{\alpha-1}} \sup_{z\geq A}\lacc\lva\Phi_{\a,\gamma,\nu_0}(z^{1/\alpha-\gamma}L_1^+) -\Phi_{\a,\gamma,\nu_0}(0+)\rva\racc\rcr\, \Pb[T_x<\infty]
\end{multline*}
and (\ref{llinf}) follows by Corollary \ref{eq:A}, the continuity of $\Phi_{\a,\gamma,\nu_0}$ at zero, and dominated convergence.

\qed

\section{Proof of Theorem B}

\subsection{The lower bound} This part is easy and relies essentially on the identity (\ref{Variance}). Introducing
$$T^{(\beta)}_x \;=\; \inf\{ t\ge 0, \, L^{(\beta)}_t = t^{\gamma+\beta} +x\}\qquad \text{and}\qquad {\widehat T}^{(\beta)}_x\;=\;\inf\{t\geq 0,\,  L^{(\beta)}_t =  (t^{\gamma +\beta} +1)\, x^{\frac{\gamma\alpha -1}{\alpha(\gamma +\beta)}} \},$$
we see by scaling that
\begin{equation}
\label{Ska}
\Pb[\Margb\ge x]\; =\; \Pb[T^{(\beta)}_x < \infty]\; =\;  \Pb[{\widehat T}^{(\beta)}_x < \infty].
\end{equation}
Setting
$$s_\ast\; =\; \arg\min\{ s^{-\beta-1/\a}(s^{\gamma+\beta}+1)\}\; =\; \left( \frac{1+\alpha \beta}{\gamma\a-1}\right)^{\frac{1}{\gamma+\beta}}$$
and
$$m_\ast\; =\; \min_{s >0}\{ s^{-\beta-1/\a}(s^{\gamma+\beta}+1), \, s >0\}\; =\;\a(\gamma+\beta) (\alpha \beta +1)^{-\frac{1+\alpha\beta}{\a(\gamma+\beta)}}(\gamma\a -1)^{\frac{1-\gamma\a}{\a(\gamma +\beta)}},$$
a further scaling argument implies
$$  \Pb[{\widehat T}^{(\beta)}_x < \infty]\;\ge\; \Pb\lcr L^{(\beta)}_{s_\ast}\geq (s_\ast^{\gamma+\beta}+1)\,  x^{\frac{\gamma\alpha -1}{\alpha(\gamma+\beta)}} \rcr\; =\; 
\Pb[L_1\geq  (1+\alpha\beta)^{1/\alpha} m_\ast\,  x^{\frac{\gamma\alpha -1}{\alpha(\gamma+\beta)}}].$$
When $L$ has positive jumps, applying e.g. Property 1.2.15 in \cite{SamTaq} yields the required lower bound
$$ \Pb[\Margb\ge x] \; \geq \; \kappa\, x^{\frac{1-\gamma\a}{\gamma+\beta}}, \qquad x\rightarrow \infty,$$
for some $\kappa > 0.$ When $L$ has no positive jumps, we obtain from (\ref{eq:iGamma}) and some simplifications the required lower bound 
$$\liminf_{x\to\infty} x^{\frac{1-\gamma\a}{(\alpha-1)(\gamma+\beta)}} \log \Pb[\M^{(\beta)}_{\a,1/\a,\gamma} \geq x]\; \ge\; - c_{\a,\beta,\gamma}.$$

\subsection{The upper bound in the case with positive jumps} Introducing the parameter
$$\delta\; =\; \frac{\gamma\alpha -1}{\alpha(\gamma +\beta)}\;\in\; (0,1)$$
and fixing $\varepsilon>0$ small enough such that $\eta = 2^{\delta}(1+\varepsilon)^{\delta-1} > 1,$ define the stochastically increasing family of stopping times
$$ {\widehat T}^{(\beta,k)}_x\;=\;\inf\{t\geq 0,\,  L^{(\beta)}_t - (1+\varepsilon)^{-k} t^{\gamma +\beta} x^\delta\, =\, 2^k x^{\delta} \},\qquad k\ge 0.$$ By (\ref{Ska}), we have the telescoping decomposition

$$\Pb[\Margb\ge x]\; =\; \Pb[{\widehat T}^{(\beta,0)}_x < \infty]\; =\; \sum_{k\ge 0}\lpa \Pb[{\widehat T}^{(\beta, k)}_x < \infty]\, -\,\Pb[{\widehat T}^{(\beta, k+1)}_x < \infty]\rpa.$$
We first consider the case $\gamma+\beta \ge 1.$ Setting $r_k=  (3 \times 2^k (1+\varepsilon)^k)^{\frac{1}{\gamma+\beta}},$ we can bound
\begin{eqnarray*}
\Pb[{\widehat T}^{(\beta, k)}_x < \infty] & \le &\Pb\lcr\sup_{t\in [0,r_k]} \{L_t^{(\beta)}\}\geq 2^k x^{\delta} \rcr\; +\; \Pb\lcr \sup_{t\geq r_k}\{ L_t^{(\beta)} - (1+\varepsilon)^{-k}t^{\gamma+\beta} x^\delta\}\geq 2^k x^{\delta}\rcr\\
& \le &  \Pb\lcr L_1^\ast \ge \eta^k  3^{\delta -1}    x^\delta   \rcr \; +\;  \Pb\lcr\sup_{t\geq 0}\{ L_{t+r_k}^{(\beta)} -(1+\varepsilon)^{-k} t^{\gamma+\beta}x^\delta\} \geq 2^{k+2}x^{\delta}\rcr, 
\end{eqnarray*}
where in the second line we have used the a.s. inequality $\sup_{t\in[0,1]}\{L_{t}^{(\beta)}\}\le L_1^\ast,$ which is obvious, and the equally obvious deterministic inequality  
\begin{equation}
\label{Dieter}
(t+r_k)^{\gamma+\beta}\; \geq\; t^{\gamma+\beta}\,+\, r_k^{\gamma+\beta}
\end{equation} 
for all $t\ge 0,$ which follows from $\gamma +\beta \ge 1.$ The next step is to write down the process decomposition
\begin{eqnarray}
\label{eq:decomp}L^{(\beta)}_{t+r_k} & = & \lpa \beta\int_0^{r_k}(t+r_k-u)^{\beta-1}\,L_u\, du \; +\; t^\beta\, L_{r_k}\rpa\;+\;  \beta\int_0^{t}(t-s)^{\beta-1}\,(L_{s+r_k}-L_{r_k})\, ds\\
\nonumber& \elaw & \lpa \beta\int_0^{r_k}(t+r_k-u)^{\beta-1}\,L_u\, du \; +\; t^\beta\, L_{r_k}\rpa\; +\; {\widehat L}^{(\beta)}_t\; \le \; (t+r_k)^\beta L^\ast_{r_k}\; +\; {\widehat L}^{(\beta)}_t
\end{eqnarray}
with $\{{\widehat L}^{(\beta)}_t, \, t\ge 0\}$ an independent copy of $\{L^{(\beta)}_t, \, t\ge 0\},$ which implies
\begin{multline*}
\Pb\lcr\sup_{t\geq 0}\{ L_{t+r_k}^{(\beta)} -(1+\varepsilon)^{-k} t^{\gamma+\beta}x^\delta\} \geq 2^{k+2}x^{\delta}\rcr
\\\!\!\!\leq\;\;  \Pb[{\widehat T}^{(\beta, k+1)}_x < \infty]\; +\; \Pb\lcr\sup_{t\geq 0} \{ L_{r_k}^\ast (t+r_k)^\beta - \varepsilon \,(1+\varepsilon)^{-k-1} t^{\gamma+\beta} x^\delta\}\geq 2^{k+1} x^{\delta}\rcr
\\\;\;\;\leq\;\;  \Pb[{\widehat T}^{(\beta, k+1)}_x < \infty]\; +\; \Pb\lcr c_\beta\, r_k^\beta L_{r_k}^\ast \,+\, \sup_{t\geq 0} \{ c_\beta\,L_{r_k}^\ast t^\beta - \varepsilon \,(1+\varepsilon)^{-k-1} t^{\gamma+\beta} x^\delta\}\geq 2^{k+1} x^{\delta}\rcr,
\end{multline*}
where $c_\beta = 2^{\vert \beta -1\vert}$ and we have used $(t+s)^\beta \le c_\beta (t^\beta + s^\beta)$ for all $t,s\ge 0.$ The second term on the right-hand side is bounded by 
\begin{multline*}
\Pb\lcr L_1^\ast \ge \eta^k  3^{\delta -1}c_\beta^{-1} x^\delta   \rcr \; +\; \Pb\lcr \sup_{t\geq 0} \{ c_\beta\,L_{r_k}^\ast t^\beta - \varepsilon \,(1+\varepsilon)^{-k-1} t^{\gamma+\beta} x^\delta\}\geq 2^k x^{\delta}\rcr
\\ =\; \Pb\lcr L_1^\ast \ge \eta^k  3^{\delta -1}c_\beta^{-1} x^\delta   \rcr\; +\; \Pb\lcr L_1^\ast \ge \eta^k \kappa\, x^\delta   \rcr
\end{multline*}
for some positive constant $\kappa$ not depending on $k,x.$ Setting ${\hat \kappa} = \min\{\kappa, 3^{\delta -1}c_\beta^{-1}\} > 0,$ and putting everything together, we finally obtain
$$\Pb[\Margb\ge x]\; \le\; 3\,\sum_{k\ge 0}\Pb\lcr L_1^\ast \ge \eta^k{\hat \kappa} \, x^\delta   \rcr\; \sim\;  \frac{3\,{\hat \kappa}^{-\a} \Gamma(\alpha)  \sin(\pi\alpha \rho)}{\pi (1-\eta^{-\a})}\, x^{\frac{1-\gamma\a}{\gamma +\beta}},$$ 
where the estimate follows at once from (\ref{Bing}) and direct summation. This completes the proof for $\gamma +\beta \ge 1.$ The case $\gamma +\beta < 1$ follows along the same lines, except that (\ref{Dieter}) is not true anymore. We hence set  
$$\lambda\; =\;\frac{\varepsilon}{2(1+\varepsilon)}\,\in\, (0,1)\qquad\mbox{and}\qquad r_k\; =\; (3\lambda^{-1} \times 2^k (1+\varepsilon)^k)^{\frac{1}{\gamma+\beta}},\quad k\ge 0.$$ 
Using the obvious inequality $(t+r_k)^{\gamma+\beta}\geq (1-\lambda)t^{\gamma+\beta}+ \lambda r_k^{\gamma+\beta}$
leads first to 
$$\Pb[{\widehat T}^{(\beta, k)}_x < \infty]\; \le\;   \Pb\lcr L_1^\ast \ge \eta^k  3^{\delta -1} x^\delta   \rcr \; +\;  \Pb\lcr\sup_{t\geq 0}\{ L_{t+r_k}^{(\beta)} -(1-\lambda)(1+\varepsilon)^{-k} t^{\gamma+\beta}x^\delta\} \geq 2^{k+2}x^{\delta}\rcr.$$
Then, we can bound
\begin{multline*}
\Pb\lcr\sup_{t\geq 0}\{ L_{t+r_k}^{(\beta)} -(1-\lambda)(1+\varepsilon)^{-k} t^{\gamma+\beta}x^\delta\} \geq 2^{k+2}x^{\delta}\rcr
\\\!\!\!\leq\;\;  \Pb[{\widehat T}^{(\beta, k+1)}_x < \infty]\; +\; \Pb\lcr\sup_{t\geq 0} \{ 2L_{r_k}^\ast (t+r_k)^\beta - \varepsilon \,(1+\varepsilon)^{-k-1} t^{\gamma+\beta} x^\delta\}\geq 2^{k+2} x^{\delta}\rcr,
\end{multline*}
and the proof is finished similarly.

\qed

\subsection{The upper bound in the case without positive jumps} The argument relies on the following well-known association lemma, which will also be used during the proof of Theorem C.

\begin{lemma}
\label{lem:ST}
Let $F, G$ be two bounded functionals on the Skorokhod space $\cD(\R^+,\R)$ being both non-increasing or both non-decreasing. Then, one has  
$$\E\lcr F(L_u,\, u\geq 0)\, G(L_u,\, u\geq 0)\rcr\; \geq\; \E\lcr F (L_u,\, u\geq 0)\rcr \E\lcr G(L_u,\, u\geq 0)\rcr.$$
\end{lemma}

\proof
By c\`ad-l\`ag approximation, it is enough to consider the case when $F,G$ depend only on a finite number of points. With the notation of Chapter 4.6 in \cite{SamTaq}, we are hence reduced to show that the random vector $(L_{t_1}, L_{t_2}, \ldots, L_{t_n})$ is associated for every $n\ge 2$ and $0 < t_1 < \ldots < t_n.$ By independence of the increments we have $(L_{t_1}, L_{t_2}, \ldots, L_{t_n}) = (X_1, X_1 +X_2, \ldots, X_1+\cdots +X_n),$ where the $X_i$'s are mutually independent real random variables, making the vector $X=(X_1, \ldots, X_n)$ trivially associated. We can then apply e.g. Exercise 4.25 p. 220 in \cite{SamTaq}.

\endproof

Let us now finish the proof. For simplicity, we will set $T_x$ for $T_x^{(\beta)}\!.$ Let $\varepsilon>0$ and fix $\delta$ small enough such that $\eta = 1- (1-\varepsilon)(\delta+1)^{\beta+\gamma} > 0.$ Using the absence of positive jumps, we obtain
\begin{align}
\notag\int_0^{\infty} \Pb[L_t^{(\beta)} \geq (1-\varepsilon)t^{\beta+\gamma} +x]\,dt & =  \int_0^{\infty} \Pb\left[L_t^{(\beta)} -L_{T_x}^{(\beta)} \geq (1-\varepsilon)t^{\beta+\gamma} - T_x^{\beta+\gamma}, \; T_x<+\infty \right]dt \\
\label{eq:nojumps1}& \ge \int_0^{\delta}  \Pb\left[L_{T_x(t+1)}^{(\beta)}-L_{T_x}^{(\beta)}  \geq -\eta\,T_x^{\beta+\gamma}, \; T_x <\infty\rcr  dt.
\end{align}
By (\ref{Variance}) and a change of variable, the left-hand side equals 
$$\int_0^{\infty} \Pb[L_t^{(\beta)} \geq (1-\varepsilon)t^{\beta+\gamma} +x]\,dt\; = \; \kappa_\varepsilon\, \int_0^{\infty} t^{-\frac{\a\beta}{1+\a\beta}}\, \Pb[L_1 \geq (1+\alpha\beta)^{1/\alpha}t^{-1/\a}(t^{\frac{\gamma+\beta}{1+\a\beta}} + c_\varepsilon x)]\,dt$$
for some positive constants $\kappa_\varepsilon, c_\varepsilon$ such that $c_\varepsilon \to 1$ as $\varepsilon\to 0$ and, by Lemma \ref{lem:Lap}, we first deduce
$$\log\int_0^{\infty} \Pb[L_t^{(\beta)} \geq (1-\varepsilon)t^{\beta+\gamma} +x]\,dt \;\sim\; - c_{\a,\beta,\gamma} (c_\varepsilon x)^{\frac{\gamma\a -1}{(\a-1)(\gamma +\beta)}}.$$
We shall now separate the proof according as $\beta\geq1$ or $\beta<1$.\\

 Assume first $\beta\geq 1$. Bounding the right-hand side of (\ref{eq:nojumps1}) leads to
$$
\int_0^{\infty} \Pb[L_t^{(\beta)} \geq (1-\varepsilon)t^{\beta+\gamma} +x]\,dt  \ge\delta\,  \Pb\left[  \inf_{u\ge 1,t\le \delta}\{ u^{-\beta-\gamma} (L_{u(t+1)}^{(\beta)}-L_u^{(\beta)} )\}\geq -\eta, \;1< T_x <\infty\rcr, 
$$
whence
\begin{multline}
\label{eq:Z}
\Pb\lcr \inf_{u\ge 1,t\le \delta}\{ u^{-\beta-\gamma} (L_{u(t+1)}^{(\beta)}-L_u^{(\beta)} )\} \geq -\eta,\; T_x <\infty\rcr\\ \le \; \delta^{-1}\int_0^{\infty} \Pb[L_t^{(\beta)} \geq (1-\varepsilon)t^{\beta+\gamma} +x]\,dt\, +\, \Pb[T_x \le 1].\qquad
\end{multline}
We next observe that the contribution of $\Pb[T_x \le 1]$ in the right-hand side of (\ref{eq:Z}) is negligible, using the obvious bound 
$$\Pb[T_x\leq 1]\;\leq\; \Pb[\tau_x\leq 1]$$ 
with $\tau_x = \inf\{t\geq0,\; L_t^{(\beta)}=x\},$ the crude estimates 
$$\log\Pb [\tau_x \le 1]\;\asymp\; \log\Pb [ L_1 > x] \;\asymp\; -x^{\frac{\alpha}{\alpha-1}}$$
and the strict inequality
$$\frac{\alpha \gamma-1}{(\alpha-1)(\beta+\gamma)} \; < \;\frac{\alpha}{\alpha-1}\cdot$$
Above, the crude estimates are a consequence of (\ref{Variance}), (\ref{eq:iGamma}) and 
$$\Pb[L_1^{(\beta)} >x ]\; \le\; \Pb[\tau_x \le 1]\; \le\; \Pb\lcr\sup_{t\le 1}\lacc L_t\racc >  x\rcr\; =\; \a\,\Pb[L_1 >  x],$$the last equality being well-known as the reflection principle for spectrally negative stable L\'evy processes - see e.g. Exercises 29.7 and 29.18 in \cite{Sato}. Finally, we notice that
$$ u^{-\beta-\gamma} (L_{u(t+1)}^{(\beta)}-L_u^{(\beta)})\; =\;\beta\int_0^{1+t} \lpa (1+t -s)^{\beta -1} - (1-s)^{\beta -1}\Un_{\{s\le 1\}}\rpa \, \frac{L_{us}}{u^\gamma}\,ds$$
is an increasing functional of $\{L_s, \, s\ge 0\},$ because $\beta \ge 1.$ Applying Lemma \ref{lem:ST}, we obtain 
\begin{multline*}
\Pb\lcr \inf_{u\ge 1,t\le \delta}\{ u^{-\beta-\gamma} (L_{u(t+1)}^{(\beta)}-L_u^{(\beta)} )\} \geq -\eta,\; T_x <\infty\rcr\\
\ge\; \Pb\lcr \inf_{u\ge 1,t\le \delta}\{ u^{-\beta-\gamma} (L_{u(t+1)}^{(\beta)}-L_u^{(\beta)} )\}\geq -\eta\rcr\;\Pb\lcr T_x <\infty\rcr\; = \; \kappa\,\Pb\lcr T_x <\infty\rcr
\end{multline*}
for some $\kappa > 0$ not depending on $x.$ Putting everything together, we get
$$\limsup_{x\to\infty} x^{\frac{1-\gamma\a}{(\a-1)(\gamma +\beta)}}\,\log \Pb[T_x <\infty]\; \le\; - c_{\a,\beta,\gamma}\, c_\varepsilon^{\frac{\gamma\a -1}{(\a-1)(\gamma +\beta)}},$$
which, letting $\varepsilon\to 0$, completes the proof in the case $\beta\geq1$.\\

Assume second $\beta< 1$. We set
$$\sigma_t= \beta  \int_0^{1} \left\{ (1-s)^{\beta-1} -(1+t -s)^{\beta-1}\right\} s^{\gamma} ds$$
which is a positive increasing function on $(0,\infty)$ such that $\sigma_t\rightarrow 0$  as $t\rightarrow 0$. 
Replacing $T_x^{\beta+\gamma}$ by
$$T_x^{\beta+\gamma} =  \frac{ \beta}{\sigma_t} \int_0^{T_x} \left\{ (T_x-s)^{\beta-1} -(T_x(1+t) -s)^{\beta-1}\right\} s^{\gamma} ds$$
we deduce using a change of variable that 
\begin{multline*}
L_{T_x(t+1)}^{(\beta)}-L_{T_x}^{(\beta)} + \frac{\eta}{2}T_x^{\beta+\gamma}\geq  \beta T_x^\beta\int_{1}^{1+t}  (1+t -s)^{\beta -1}L_{sT_x} ds
-  T_x^\beta h_\beta(t) \sup_{u\geq0}\left\{L_u- \frac{\eta}{2\sigma_t}u^\gamma\right\}.
\end{multline*}
where $h_\beta(t)=1+t^\beta -(1+t)^\beta$ is increasing in $t$.
Going back to (\ref{eq:nojumps1}), and taking $a<\delta$, the  right-hand side is then greater than :
\begin{equation}\label{eq:a}
a \Pb\left[ F_\delta(L)-  \frac{ h_\beta(\delta)}{ \delta^{\gamma}x^{\frac{\gamma}{\beta+\gamma}}} \sup_{s\geq0} \left\{L_{s} - \frac{\eta}{2 \sigma_a} s^\gamma\right\}  \geq -\eta/2, \;\delta x^{\frac{1}{\beta+\gamma}} <T_x <\infty\rcr 
\end{equation}
where
$$
F_\delta(L)= \beta   \inf_{t\leq \delta}  \int_1^{1+t}  (1+t -s)^{\beta -1}  \inf_{u\geq 1}\frac{L_{su}}{u^\gamma}\,ds
$$
is an increasing functional of $L$. We next observe that, cutting (\ref{eq:a}) in two as in (\ref{eq:Z}), the second term will be negligible by taking $\delta$ small enough since
$$\Pb[T_x\leq  \delta x^{\frac{1}{\gamma+\beta}}]\;\leq\; \Pb[\tau_x\leq  \delta  x^{\frac{1}{\gamma+\beta}}]$$ 
and 
$$\log\Pb [\tau_x \le \delta x^{\frac{1}{\gamma+\beta}}]\;\asymp\; \log\Pb [ L_1 > \delta^{-\beta - \frac{1}{\alpha}}   x^{\frac{\alpha\gamma-1}{\alpha(\gamma+\beta)}}] \;\asymp\; - \delta^{- \frac{\alpha\beta+1}{\alpha-1}} x^{\frac{\alpha\gamma-1}{(\alpha-1)(\gamma+\beta)}}.$$
Thus, it remains to deal with the term :
$$   \Pb\left[ F_\delta(L) \geq -\eta/4,\, T_x <\infty\rcr  -  \Pb\left[  h_\beta(\delta)  \sup_{s\geq0} \left\{L_{s} - \frac{\eta}{2 \sigma_a} s^\gamma\right\} \geq \frac{\eta}{4} \delta^{\gamma}x^{\frac{\gamma}{\beta+\gamma}} \right]. 
$$
From Theorem A and using the scaling of $L$, the second term behaves as
$$  \log \Pb\left[ h_\beta(\delta) \left(\sigma_a\right)^{\frac{1}{\alpha\gamma-1}}  \sup_{s\geq0} \left\{L_{s} -\frac{\eta}{2} s^\gamma\right\} \geq\frac{\eta}{4} \delta^\gamma  x^{\frac{\gamma}{\beta+\gamma}} \right]\asymp - \left(\sigma_a\right)^{-\frac{1}{\gamma(\alpha-1)}}  x^{\frac{\gamma \alpha-1}{(\alpha-1)(\gamma+\beta)}}$$
which is negligible by taking $a$ small enough. The proof is then concluded as in the case $\beta\geq1$ by applying Lemma \ref{lem:ST} to the term $
\Pb\left[ F_\delta(L) \geq -\eta/4, \; T_x <\infty\rcr$.\\
\qed

\begin{remark} In the particular case $\beta =1$ of the integrated stable process, we may proceed as in the proof of Theorem A, and obtain a more precise upper bound. The strong Markov property at $T_x$ for the two-dimensional process 
$$\{(L^{(1)}_t, L_t), \, t\ge 0\},$$ 
a scaling argument and (\ref{Variance}) imply firstly
\begin{multline*}
\int_0^\infty t^{\nu_0}\, \Pb[L^{(1)}_t> t^{1+\gamma} +x]\, dt \\
 = \;\E\left[\Un_{\{T_x<\infty\}} \, T_x^{1+\nu_0}\, \int_0^{\infty} (1+t)^{\nu_0}\, \Un_{\{\tLL_1^+ + (tT_x)^{-1/\a}  (L_{T_x} -(1+\gamma) T_x^{\gamma}) >  T_x^{\gamma-1/\alpha}\,\psi_{\a,\gamma} (t)\}}\, dt\rcr
\end{multline*}
where $\psi_{\a,\gamma}(t)= t^{-1-1/\a}((t+1)^{1+\gamma} -1 - (1+\gamma) t)$ is again an homeomorphism from $(0,\infty)$ to $(0,\infty),$ and $T_x$ and $\tLL_1^+$ are independent. We can then bound
$$
\int_0^\infty t^{\nu_0}\, \Pb[L^{(1)}_t> t^{1+\gamma} +x]\, dt\;\geq \;\E\left[\Un_{\{T_x<\infty\}} \, T_x^{1+\nu_0}\, \int_0^{\infty} (1+t)^{\nu_0}\, \Un_{\{\tLL_1^+ >  T_x^{\gamma-1/\alpha}\,\psi_{\a,\gamma} (t)\}}\, dt\rcr,$$
using the crucial fact that the derivative of $t\mapsto L^{(1)}_t - t^{1+\gamma}$ at $T_x,$ which equals $L_{T_x} - (1+\gamma) T_x^{\gamma},$ is a.s. non-negative. This leads to
$$ \int_0^{\infty} t^{\nu_0}\,\Pb[L_1^+> t^{-1-1/\a}(t^{1+\gamma} +x)]\,dt\; \ge\;\frac{1}{1+\nu_0}\; \E\lcr\Un_{\{T_x<\infty\}}\left(\tLL_1^+\right)^{\frac{(1+\nu_0)\alpha}{\gamma\a-1}} \Psi_{\a,\gamma} (T_x^{1/\a-\gamma}\,\tLL_1^+)\rcr$$
where 
$$\Psi_{\a,\gamma}(x)\; =\; x^{-\frac{(1+\nu_0)\alpha}{\gamma\a-1}} \left((1+\psi_{\a,\gamma}^{-1}(x))^{1+\nu_0}-1\right)$$
is again bounded away from zero and $\infty,$ by the fateful choice of $\nu_0.$ We finally obtain 
$$\Pb[T_x <\infty]\; \le \; \kappa_+ \, \int_0^{\infty} t^{\nu_0}\,\Pb[L_1^+> t^{-1-1/\a}(t^{1+\gamma} +x)]\,dt$$
for some $\kappa_+ \in(0,\infty),$ and an appropriate modification of Lemma 1 yields  
$$\Pb[\M^{(1)}_{\a,1/\a,\gamma} \geq x]\; \le\; {\hat \kappa}_+\,\exp\lacc- c_{\a,1,\gamma}\,x^{\frac{\gamma\a-1}{(\alpha-1)(1+\gamma)}}\racc$$
at infinity, for some other ${\hat \kappa}_+ \in(0,\infty).$ Unfortunately, the precise lower bound which can be derived from (\ref{eq:iGamma}) is different: one gets
$$\Pb[\M^{(1)}_{\a,1/\a,\gamma} \geq x]\; \ge\; {\hat \kappa}_-\,x^{\frac{1-\gamma\a}{2(\alpha-1)(1+\gamma)}}\,\exp\lacc- c_{\a,1,\gamma}\,x^{\frac{\gamma\a-1}{(\alpha-1)(1+\gamma)}}\racc$$
for some ${\hat \kappa}_-\in(0,\infty),$ and the exact polynomial speed before the exponential term remains unknown. We believe that this speed is given in the lower bound, and we refer to Remark \ref{Precise} (c) below for a general conjecture.

\end{remark}

\subsection{A more precise estimate in the Brownian case} In this paragraph we specify the general results of \cite{HP} to the process $L^{(\beta)}$ in the case $\alpha =2,$ and we get a refinement of Theorem B (b). Observe that in this framework we can also consider the wider range $\beta > -1/2.$ It turns out that a transition phenomenon occurs around $\beta =1/2.$

\begin{proposition}
\label{ApplyHP}
Assume $\gamma>1/2.$

\medskip 

{\em (a)} If $\beta \in(-1/2, 1/2),$ there exists $\kappa_{\beta,\gamma} > 0$ such that
$$\Pb[\M_{2,1/2,\gamma}^{(\beta)} \geq x] \;\sim\; \kappa_{\beta,\gamma}\, x^{\frac{2\beta(1-2\gamma)}{(2\beta +1)(\gamma+\beta)}}\, \exp\lacc -c_{2,\beta,\gamma}\, x^{\frac{2\gamma -1}{\gamma+\beta}}\racc.$$

{\em (b)} If $\beta > 1/2,$ there exists ${\tilde \kappa_{\beta,\gamma}} > 0$ such that
$$\Pb[\M_{2,1/2,\gamma}^{(\beta)} \geq x] \;\sim\; {\tilde \kappa_{\beta,\gamma}}\, x^{\frac{1-2\gamma}{2(\gamma+\beta)}}\, \exp\lacc -c_{2,\beta,\gamma}\, x^{\frac{2\gamma -1}{\gamma+\beta}}\racc.$$

\end{proposition}
 
\proof

With our normalization, one has $L_1 \sim \cN(0, 2)$ and a scaling argument implies
$$\M_{2,1/2,\gamma}^{(\beta)}\;\elaw\; \left(\frac{2}{2\beta+1}\right)^{\frac{\gamma +\beta}{2\gamma-1}}\, {\widetilde \M_{\beta,\gamma}}$$
where
$${\widetilde \M_{\beta,\gamma}}\; =\; \sup_{t>0}\lacc \sqrt{2\beta+1}\int_0^t (t-s)^\beta \, dB_s\, - \, t^{\gamma +\beta}\racc$$
and $\{B_t,\, t\ge 0\}$ is a standard linear Brownian motion. Setting $H = \beta +1/2,$ the process
$$X_t\;=\;  \sqrt{2\beta +1}\int_0^t (t-s)^\beta \, dB_s, \qquad t \ge 0,$$
is Gaussian with mean $0$ and variance $t^{2H},$ and self-similar with index $H.$ With the notation of Section 1 in \cite{HP}, we have
\begin{equation}
\label{Een}
s_0\; =\; \lpa\frac{2\beta +1}{2\gamma -1}\rpa^{\frac{1}{\gamma +\beta}}\qquad\mbox{and}\qquad A\; =\; 2 \sqrt{\frac{c_{2,\beta,\gamma}}{2\beta+1}}.
\end{equation}
We now wish to apply Theorem 1 in \cite{HP}, whose statement deals with case $H\in (0,1)$ but we can actually consider any $H >0$ by Remark 1 therein. Following (7) in \cite{HP}, the next step is to evaluate the behaviour of $\E[(Y_t - Y_s)^2]$ as $t,s\to s_0,$ having set $Y_t = t^{-H}X_t$ for all $t > 0.$ Because of the time normalization, we have not found any precise reference for this behaviour in the literature and so we give the details. For $0 < s <t,$ we compute
$$\E[(Y_t - Y_s)^2]\; =\; 2 \; -\; I_\beta (x)$$
with $x=st^{-1}\in (0,1)$ and
$$I_\beta(x)\; =\; 2(2\beta +1)\sqrt{x}\int_0^1 (1-u)^\beta (1-xu)^\beta\, du.$$
We need to study the asymptotic behaviour of $I_\beta(x)$ as $y = 1-x\to 0.$ If $\beta > 1/2,$ rewriting
$$I_\beta(x)\; =\; 2(2\beta + 1)\sqrt{1-y}\int_0^1 (1-u)^{2\beta} (1+yu(1-u)^{-1})^\beta\, du,$$
making a Taylor expansion of order 2 of both quantities in $y$ and evaluating the two underlying Beta integrals leads to
$$I_\beta(x)\; =\; 2\; -\; \frac{(2\beta +1) y^2}{4(2\beta -1)}\; +\; o(y^2).$$
This shows that (7) in \cite{HP} holds with
\begin{equation}
\label{Twee}
\alpha\; = \;2 \qquad \mbox{and}\qquad D\; =\; \frac{(2\beta +1)}{4(2\beta -1)s_0^2}\cdot
\end{equation}
If $\beta < 1/2,$ the argument does not apply because the second Beta integral diverges. We first rewrite
$$I_\beta(x)\; =\; \frac{2(2\beta+1)\sqrt{x}}{\beta +1}\;\pFq{2}{1}{-\beta,1}{\beta +2}{x}\; =\;\frac{2(2\beta+1)\sqrt{x} y^\beta}{\beta +1}\;\pFq{2}{1}{-\beta,\beta + 1}{\beta +2}{-xy^{-1}},$$
where the first equality follows from Euler's formula and the second one from Pfaff's transformation for the hypergeometric function - see respectively 2.1.3(10) and 2.1.4(22) in \cite{EMOT}. Applying next the residue transformation 2.1.4(17) in \cite{EMOT}, we obtain
\begin{eqnarray*}
I_\beta(x) & = & 2 x^{\beta +1/2}\pFq{2}{1}{-\beta,-1-2\beta}{-2\beta}{-yx^{-1}}\; -\; \frac{2\Gamma(\beta +1)\Gamma(-2\beta)}{\Gamma(-\beta)} \, x^{-\beta -1/2} y^{2\beta +1}\\
& = &  2 \; -\; \frac{\Gamma^2(\beta +1)}{\Gamma(2\beta +1)\,\cos(\pi\beta)} \, y^{2\beta +1} \; +\; O(y^2).
\end{eqnarray*}
This shows that (7) in \cite{HP} holds with
\begin{equation}
\label{Drii}
\alpha\; = \;2\beta +1  \qquad \mbox{and}\qquad D\; =\; \frac{\Gamma^2(\beta +1)}{\Gamma(2\beta +1)\,\cos(\pi\beta)\, s_0^{2\beta +1}}\cdot
\end{equation}
Putting (\ref{Een}) and (\ref{Twee}) resp. (\ref{Drii}) together with (10) resp. (9) in \cite{HP}, using the standard estimate $\sqrt{2\pi}\Psi (u)\sim u^{-1} e^{-u^2/2}$ for the tail of the unit normal distribution, and proceeding to the necessary simplifications, we obtain our required asymptotics with the two different regimes.
  
\endproof

\begin{remark} 
\label{Precise}
(a) For $\beta= 1/2,$ the transformation 2.1.4(18) in \cite{EMOT} with $m=2$ exhibits a logarithmic term: one has the non-trivial closed formula
\begin{eqnarray*}
I_{1/2}(x) & =& 2\; +\;\frac{y^2}{2(1-y)} (\psi(3/2) \,-\,\psi(3)\, +\,\log(y)\,-\,\log(1-y))
\end{eqnarray*}
where $\psi$ is the digamma function. This implies
$$\E[(Y_t-Y_s)^2]\; \sim\; -\frac{(t-s)^2\log\vert t-s\vert}{2 s_0^2}$$
as $t,s\to s_0,$ and we cannot apply the results of \cite{HP}. We believe that 
$$\Pb[\M_{2,1/2,\gamma}^{(1/2)} \geq x] \;\sim\;\kappa\, (\log x)^\delta \,x^{\frac{1-2\gamma}{2\gamma+1}}\, \exp\lacc -c_{2,1/2,\gamma}\, x^{\frac{2\gamma -1}{(\gamma+1/2)}}\racc$$
for some $\kappa > 0$ and $\delta\neq 0$ to be determined, the logarithmic correction being heuristically due to the 1-self-similarity of 
$$t\;\mapsto\;\int_0^t \sqrt{t-s}\; dB_s.$$

(b) The constants $\kappa_{\beta,\gamma}$ and ${\tilde \kappa}_{\beta,\gamma}$ can also be evaluated from Theorem 1 in \cite{HP}, but they have a complicated form in general. For $\beta > 1/2,$ one gets
$${\tilde \kappa}_{\beta,\gamma}\; =\; \sqrt{\frac{2(2\beta\gamma -\beta-\gamma +1)}{\pi\,c_{2,\beta,\gamma}\,(2\beta-1)(2\gamma -1)}}\cdot$$
For $\beta\in(-1/2,1/2),$ one obtains
$$\frac{H_{2\beta +1}}{\sqrt{(2\beta +1)(2\gamma-1)}} (\gamma+\beta)^{-\frac{4\beta}{2\beta +1}}  \left(\frac{2\gamma-1}{2\beta+1}\right)^{\frac{2\beta(2\gamma-1)}{(2\beta+1)(\gamma+\beta)}} \lpa \frac{\Gamma^2(\beta +1)}{\Gamma(2\beta +1)\,\cos(\pi\beta)}\rpa^{\frac{1}{2\beta +1}}$$
where $\{B_H(t),\, t\ge 0\}$ is a standard fractional Brownian motion with Hurst parameter $H=\beta +1/2,$ and
$$H_{2\beta +1}\; =\; \lim_{T\to \infty} \frac{1}{T}\,\E\lcr \exp\lacc \max_{0\le t\le T} (\sqrt{2} B_H(t) - t^{2H})\racc\rcr$$
is, in the words of \cite{HP}, a ``well-known constant''. It does not seem to the authors that the latter constant is explicit, save for $\beta =0$ where the reflection principle and Laplace's method yield
$$\E\lcr \exp\lacc \max_{0\le t\le T} (\sqrt{2} B_t - t)\racc\rcr\; =\; 1\, +\,\frac{T^{3/2}}{2\sqrt{\pi}}\int_0^1 \sqrt{s}\lpa\int_0^\infty x\, e^{-\frac{sT (x-1)^2}{4}}\, dx\rpa ds\;\sim\; T,$$
so that $H_1 = 1$ and 
$$\kappa_{0,\gamma}\; =\; \frac{1}{\sqrt{2\gamma -1}}$$
in accordance with Theorem A (b).\\
  
(c) From Proposition \ref{ApplyHP}, it is plausible to conjecture that for $\a\in(1,2)$ one has
$$\Pb[\M^{(\beta)}_{\a,1/\a,\gamma} \geq x]\; \sim\; \kappa_{\a,\beta,\gamma}\,x^{\frac{\a\beta(1-\gamma\a)}{(\a-1+\a\beta)(\alpha-1)(\gamma +\beta)}}\,\exp\lacc- c_{\a,\beta,\gamma}\,x^{\frac{\gamma\a-1}{(\alpha-1)(\gamma +\beta)}}\racc$$
if $\beta \in (0,1-1/\a)$ and 
$$\Pb[\M^{(\beta)}_{\a,1/\a,\gamma} \geq x]\; \sim\; {\tilde \kappa_{\a,\beta,\gamma}}\,x^{\frac{1-\gamma\a}{2(\alpha-1)(\gamma +\beta)}}\,\exp\lacc- c_{\a,\beta,\gamma}\,x^{\frac{\gamma\a-1}{(\alpha-1)(\gamma +\beta)}}\racc$$
if $\beta > 1 -1/\a,$ where $\kappa_{\a,\beta,\gamma}$ and ${\tilde \kappa_{\a,\beta,\gamma}}$ are some positive and finite constants.
 
\end{remark}

\section{Proof of Theorem C}

Following the notation of \cite{PSInt}, we will set 
$$\theta\; =\; \frac{\rho}{\a(1-\rho) +1}$$
once and for all. The upper bound follows easily from 
$$\Pb\lcr\sup_{0\leq t\leq 1}\lacc L_t^{(1)} + \mu t^{1+\gamma}\racc\leq\, \varepsilon\rcr\;\leq \; \Pb\lcr\sup_{0\leq t\leq 1}\lacc L_t^{(1)}\racc \leq \varepsilon\rcr\; \leq\; \kappa\, \varepsilon^{\frac{\alpha\theta}{\alpha+1}}$$
for some $\kappa\in (0,\infty),$ where the first inequality follows from $\mu\ge 0$ and the second one from Theorem A in \cite{PSInt} and scaling. \\

The lower bound is more involved and we will need the strong Markovian character of the two-dimensional process $\{(L^{(1)}_t\!,L_t), \, t\ge 0\},$ setting by $\Pb_{(x,y)}$ for its law starting from $(x,y)\in \R^2.$ Define the stopping time 
$$R_{\varepsilon}\; =\; \inf\left\{t\geq0,\, L_t^{(1)} + \mu \varepsilon^{\frac{\gamma\alpha -1}{\alpha+1}} t^{1+\gamma} =0\right\}$$
and observe first that, by scaling and translation,
$$\Pb\lcr\sup_{0\leq t\leq 1}\lacc L_t^{(1)} + \mu t^{1+\gamma}\racc\le\,\varepsilon\rcr\; =\; \Pb_{(-1,0)}\lcr R_\varepsilon\,\geq\,\varepsilon^{-\frac{\alpha}{\alpha+1}}\rcr.$$
Notice also that $\Pb_{(x,y)}\lcr R_\varepsilon \leq R_0\rcr=1$ for every $x < 0$ and $y\in\R,$ because $\mu \ge 0.$ Applying the strong Markov property at $R_\varepsilon,$ we obtain
$$\Pb_{(-1,0)}\lcr R_0 \geq 2\varepsilon^{-\frac{\alpha}{\alpha+1}}\rcr \; = \;
\E_{(-1,0)}\lcr\Pb_{(-\mu\varepsilon^{\frac{\gamma\alpha -1}{\alpha+1}} R_\varepsilon^{1+\gamma}, L_{R_\varepsilon})}\lcr R_0+ x\geq 2\varepsilon^{-\frac{\alpha}{\alpha+1}}\rcr_{\{x=R_\varepsilon\}}\rcr
$$
whose right-hand side is, by comparison, smaller than 
$$\Pb_{(-1,0)}\lcr R_\varepsilon \geq \varepsilon^{-\frac{\alpha}{\alpha+1}}\rcr \; +\; 
\E_{(-1,0)}\lcr\Un_{\{ R_\varepsilon \leq \varepsilon^{-\frac{\alpha}{\alpha+1}}\}}\,\Pb_{(-\mu\varepsilon^{\frac{\gamma\alpha -1}{\alpha+1}} R_\varepsilon^{1+\gamma}, - \mu(1+\gamma)\varepsilon^{\frac{\gamma\alpha -1}{\alpha+1} }R_{\varepsilon}^{\gamma})}\lcr R_0\geq \varepsilon^{-\frac{\alpha}{\alpha+1}}\rcr\rcr.$$
Indeed, the derivative of $t\mapsto L^{(1)}_t+ \mu \varepsilon^{\frac{\gamma\alpha -1}{\alpha+1}} t^{1+\gamma}$ at $R_\varepsilon$ equals $L_{R_\varepsilon} + \mu(1+\gamma)\varepsilon^{\frac{\gamma\alpha -1}{\alpha+1} }R_{\varepsilon}^{\gamma}$ and is a.s. non-negative under $\Pb_{(-1,0)}$. On the other hand, a further scaling argument shows that
$$\Pb_{(- x,-y)}[R_0\geq t] \;=\; \Pb_{(- 1, -yx^{-\frac{1}{\alpha+1}})}\lcr x^{\frac{\alpha}{\alpha+1}}  R_0\geq t  \rcr$$
for every $x,y,t \ge 0.$ If we now assume $\mu\le 1$ this implies, again by comparison,
\begin{multline*}
\E_{(-1,0)}\lcr\Un_{\{ R_\varepsilon \leq \varepsilon^{-\frac{\alpha}{\alpha+1}}\}}\,\Pb_{(-\mu\varepsilon^{\frac{\gamma\alpha -1}{\alpha+1}} R_\varepsilon^{1+\gamma}, - \mu(1+\gamma)\varepsilon^{\frac{\gamma\alpha -1}{\alpha+1} }R_{\varepsilon}^{\gamma})}\lcr R_0\geq \varepsilon^{-\frac{\alpha}{\alpha+1}}\rcr\rcr\\
=\;\E_{(-1,0)}\lcr\Un_{\{ R_\varepsilon \leq \varepsilon^{-\frac{\alpha}{\alpha+1}}\}}\,\Pb_{(-1,  - \mu^{\frac{\alpha}{\alpha+1}} (1+\gamma)(\varepsilon^{\frac{\alpha}{\alpha+1}} R_\varepsilon) ^{\frac{\gamma\alpha -1}{\alpha+1}})}\lcr \left(\mu x^{1+\gamma}\,\varepsilon^{\frac{\gamma\alpha -1}{\alpha+1}}\right)^{\frac{\alpha}{\alpha+1}} R_0\geq \varepsilon^{-\frac{\alpha}{\alpha+1}}\rcr_{\{x=R_\varepsilon\}}\rcr\\
 \leq \; \E_{(-1,0)}\lcr\Pb_{(-1,  -1-\gamma)}\lcr x^{\frac{\alpha(\gamma+1)}{\alpha+1}} R_0\,\geq \mu^{-\frac{\alpha}{\alpha+1}}  \varepsilon^{-\frac{\a^2(\gamma+1)}{(\alpha+1)^2}}\rcr_{\{x=R_\varepsilon\}}\rcr\\
\le\; \Pb_{(-1,  -1-\gamma)}\lcr {\widehat R_0}^{\frac{\alpha(\gamma+1)}{\alpha+1}} R_0\,\geq \mu^{-\frac{\alpha}{\alpha+1}}  \varepsilon^{-\frac{\alpha^2(\gamma +1)}{\alpha+1}}\rcr
\end{multline*}
where ${\hat R_0}$ is an independent copy of $R_0.$ Putting everything together, we obtain
$$\Pb_{(-1,0)}\lcr R_\varepsilon \geq \varepsilon^{-\frac{\alpha}{\alpha+1}}\rcr\; \ge\;\Pb_{(-1,0)}\lcr R_0 \geq 2\varepsilon^{-\frac{\alpha}{\alpha+1}}\rcr \; -\;  \Pb_{(-1,  -1-\gamma)}\lcr {\widehat R_0}^{\frac{\alpha(\gamma+1)}{\alpha+1}} R_0\,\geq \mu^{-\frac{\alpha}{\alpha+1}}  \varepsilon^{-\frac{\alpha^2(\gamma +1)}{\alpha+1}}\rcr.$$
Now by Theorem A in \cite{PSInt} we have 
$$\Pb_{(-1,  x)}\lcr R_0 > t\rcr \;\asymp\; t^{-\theta}$$
for every $x\in\R$ and since  $\alpha(\gamma+1) > \alpha+1,$ we can also infer from  Lemma 2 in \cite{PSWind} that 
$$\Pb_{(-1,  -1-\gamma)}\lcr {\widehat R_0}^{\frac{\alpha(\gamma+1)}{\alpha+1}} R_0\,\geq \, t \rcr\;\asymp\; t^{\frac{\theta(\alpha+1)}{\alpha(\gamma+1)}}.$$
This implies that there exists two finite constants $\kappa_2 \ge \kappa_1 > 0$ independent of $\mu,\varepsilon$ such that
$$\Pb_{(-1,0)}\lcr R_\varepsilon \geq \varepsilon^{-\frac{\alpha}{\alpha+1}}\rcr\; \ge\;\kappa_1\,\varepsilon^{\frac{\alpha \theta}{\alpha+1}}\, -\,  \kappa_2\, \mu^{\frac{\theta}{\gamma+1}}  \varepsilon^{\frac{\alpha \theta}{\alpha+1}},$$
which completes the proof of the lower bound for $\mu \leq\mu_0$ with $\mu_0 =  (\kappa_1/2\kappa_2)^{(\gamma+1)/\theta}> 0.$ 

Assuming finally $\mu > \mu_0$ and setting ${\bar \mu}=\mu^{\frac{\alpha}{\alpha\gamma-1}}$ and  ${\bar \mu_0}=\mu_0^{\frac{\alpha}{\alpha\gamma-1}}$ for simplicity, we have
\begin{eqnarray*}
\Pb\lcr\sup_{0\leq t\leq 1}\lacc L_t^{(1)} + \mu t^{1+\gamma}\racc\le\,\varepsilon\rcr &= &\Pb\lcr\sup_{0\leq t\leq {\bar \mu}}\lacc L_t^{(1)} +  t^{1+\gamma}\racc\le\,\varepsilon {\bar \mu}^{\frac{\alpha+1}{\alpha}} \rcr \\
&\geq &\Pb\lcr\sup_{0\leq t\leq {\bar \mu}}\lacc L_t^{(1)} +  t^{1+\gamma}\racc\le\,\varepsilon {\bar \mu_0}^{\frac{\alpha+1}{\alpha}} \rcr\\
& \geq & \Pb\lcr \sup_{{\bar \mu_0}\leq t\leq {\bar \mu}}\lacc L_t^{(1)} + t^{1+\gamma}\racc\le\,0, \;\sup_{0\leq t\leq {\bar \mu_0}}\lacc L_t^{(1)} + t^{1+\gamma}\racc\le\,\varepsilon {\bar \mu_0}^{\frac{\alpha+1}{\alpha}}\rcr\\
& \geq &\Pb \lcr \sup_{{\bar \mu_0}\leq t\leq {\bar \mu}}\lacc L_t^{(1)} + t^{1+\gamma}\racc\le\,0 \rcr \Pb\lcr\sup_{0\leq t\leq {\bar \mu_0}}\lacc L_t^{(1)} + t^{1+\gamma}\racc\le\,\varepsilon {\bar \mu_0}^{\frac{\alpha+1}{\alpha}}\rcr\\
& = & \Pb \lcr \sup_{{\bar \mu_0}\leq t\leq {\bar \mu}}\lacc L_t^{(1)} + t^{1+\gamma}\racc\le\,0 \rcr \Pb\lcr\sup_{0\leq t\leq 1}\lacc L_t^{(1)} + \mu_0 t^{1+\gamma}\racc\le\,\varepsilon\rcr\\
& \ge & \kappa\,  \varepsilon^{\frac{\alpha \theta}{\alpha+1}}
\end{eqnarray*}
for some $\kappa > 0,$ where the first and fifth equalities are obtained by scaling, the fourth inequality follows from Lemma \ref{lem:ST}, and the last inequality is a consequence of the strict positivity of ${\bar \mu_0}.$ This completes the proof.

\qed

\begin{remark} Using the same argument and Lemma VIII.4 in \cite{Ber}, one can show the following lower tail probabilities estimate for the L\'evy stable process with a positive power drift. If $\a\gamma > 1$ and $\rho\in (0,1),$ then for every $\mu \ge 0$ one has
\begin{equation}
\label{SmallV}
\Pb\lcr\sup_{0\leq t\leq 1}\lacc L_t + \mu t^{\gamma}\racc\leq\, \varepsilon\rcr\;\asymp\; \varepsilon^{\alpha\rho}.
\end{equation}
We leave the detail, which is simpler than the above, to the interested reader. This estimate for small values echoes the persistence result for large time obtained in Theorem 1 of \cite{AK}, which reads
\begin{equation}
\label{Persi}
\Pb\lcr\sup_{0\leq t\leq T}\lacc L_t + \mu t^{\gamma}\racc\leq\, 1\rcr\;=\; T^{-\rho +o(1)}
\end{equation}
for every $\mu \ge 0,$ with $\a\gamma < 1, \rho\in (0,1),$ and under the additional assumption $\a\in (0,1).$ Observe that in the absence of self-similarity, the estimates (\ref{SmallV}) and (\ref{Persi}) are different ones and cannot be deduced from one another, save for $\mu =0.$ 
\end{remark}

\bigskip

\noindent
{\bf Acknowledgement.} This work was initiated during a stay at Academia Sinica, Taipei, of the second author who would like to thank Chii-Ruey Hwang and Ju-Yi Yen for their hospitality, as well as Chien-Hao Huang and Zakhar Kabluchko for some useful discussions.

\end{document}